\documentclass[11pt,twoside]{article}

\usepackage{amsfonts}
\usepackage{amssymb}
\usepackage{amsmath}
\usepackage{graphicx}
\usepackage{hyperref}
\numberwithin{equation}{section}
\newtheorem{Prop}{\bf Proposition}[section]

\newtheorem{defn}{\bf Definition}[section]
\newtheorem{Rem}{\bf Remark}[section]
\newtheorem{Ex}{\bf Example}[section]
\newtheorem{Th}{Theorem}[section]

\begin{document}
\def \b{\Box}
\def \to{\mapsto}
\def \e{\varepsilon}
\begin{center}
{\Large {\bf  On three generalizations of the group concept: groupoid, generalized group and almost  groupoid }}
\end{center}

\begin{center}
{\bf Gheorghe IVAN}
\end{center}

\setcounter{page}{1}
\pagestyle{myheadings}

{\small {\bf Abstract}. The aim of this paper is to describe the definitions and main properties of three generalizations of the group concept, namely: groupoid, generalized group and almost groupoid. Some constructions of these algebraic structures and 
corresponding examples are presented..
{}

\section{Introduction}
\indent\indent The concept of groupoid was first introduced by H. Brandt (\cite{brand}) in 1926 and it is developed by P.J. Higgins in \cite{higgi71}.From this setting a groupoid (in the sense of Brandt)
can be thought of as a generalization of a group in which only certain multiplications are possible and it contains several  units elements. For basic results of Brandt groupoids, see (\cite{giv1b,  mack87, mapi22, rare01, west71}).\\
The concept of generalized group was introduced by M.R. Molaei  ((\cite{mola98})) in 1998.

After the introduction of topological and differentiable groupoids by Ehresmann((\cite{ehre50})  in 1950, they have been studied by many mathematicians with different approaches (\cite{brown06}). 
One of these different approaches is the structured groupoid, which is obtained by adding another algebraic structure such that the composition of groupoid is compatible with the operation of the added algebraic structure. 
A famous example of structured groupoid is the concept of a group-groupoid, defined by R. Brown and C. B. Spencer (\cite{brspen}). 

Another important objective of studying groupoids and generalized groups has been to add topological or differentiable structures, obtaining structured groupoids or structured generalized groups. 
Among the structured groupoids studied in various works, we mention:\\
Also, among the structured generalized groups studied in various works, we mention\\

The third generalization of the group concept is that of the almost groupoid. Recently, the almost groupoid concept was introduced by Mihai Ivan in 2023 (\cite{miva23a}).   It is very important
 to note that the groups category  is a subcategory of the category of almost groupoids and the almost groupoids category is a subcategory of groupoids category.

The paper is organized as follows. The Section 2 is devoted to giving fundamental definitions and basic concepts related to Brandt  groupoids.The definition of generalized group and some of its basic properties are given in Section 3. 
Section 4 present: the definitions of almost  groupoid, almost groupoid morphism and algebraic substructures of a given almost groupoid. In the final part, three constructions of new almost groupoids are presented, starting from two almost groupoids. 
For the construction of the semidirect product of two almost groupoids, Theorem 4.1 is proved. \\[-0.7cm]

\section{ Groupoids: definitions and basic results}

This section  is devoted to giving basic definitions and some results related to Brandt's groupoids. 
{\bf 2.1.~The concept  of groupoid as algebraic structure}
\begin{defn}
 (\cite{wein96})~ {\rm A  (Brandt)  \textit{groupoid $ \Gamma~$  over  $\Gamma_{0}~$}  or a {\it $ \Gamma_{0}-$groupoid}  is a pair $ (\Gamma, \Gamma _0) $ of nonempty sets such that $
\Gamma _0\subseteq \Gamma  $ endowed with the surjective maps  $~\alpha :\Gamma \rightarrow~ \Gamma _0, ~x \rightarrow ~\alpha(x)~$  and $~\beta : \Gamma~\rightarrow~\Gamma _0,~x\rightarrow \beta(x),$
$~ \mu:\Gamma_{(2)}\rightarrow \Gamma ,~(x,y)\longmapsto  \mu ( x,y):=x\cdot y,~$  where\\
$~\Gamma_{(2)}:= \{ (x,y)\in \Gamma\times \Gamma | \beta(x) = \alpha(y) \}~$ and an injective map  $~\iota:\Gamma \to \Gamma,~x \to \iota (x): = x^{ - 1},~$satisfying the following properties:\\
$(1)~$({\it associativity})~$ (x y) z = x (y z),~$ in the sense that, if one side of the equation is defined so is the other one and then they are equal;\\
$(2)~$({\it identities}): for each $ x\in \Gamma~$ follows  $~(\alpha (x),x),~(x,\beta (x))\in \Gamma_{(2)}~$ and  we have $~ \alpha(x) x = x\beta(x) = x;$\\
$(3)~$({\it inverses}): for all $~x \in \Gamma~$ follows $~(x ,x^{-1}),~(x^{-1}, x) \in \Gamma_{(2)} $ and the following equalities hold $~x \cdot x^{- 1} = \alpha (x)~$ and $~x^{-1} \cdot x = \beta(x).$}
\end{defn}
The elements of $~\Gamma_{(2)}~$ are called {\it composable pairs} of $~\Gamma.~$ The element $~\alpha(x)~$ (resp. $~\beta(x)~$) is the {\it left unit}  (resp. {\it right unit})
of $~x\in \Gamma.~$ The subset $~\Gamma_{0} = \alpha(\Gamma) = \beta(\Gamma)~$ of $~\Gamma~$  is called the {\it unit set} of $~\Gamma.~$

 A groupoid $ \Gamma $ over $~\Gamma_{0} $ with the \textit{structure functions} $ \alpha $ ({\it source}), $ \beta $ ({\it target}), $ \mu $ ({\it
partial multiplication}) and $\iota$ ({\it inversion}), is denoted by $ (\Gamma, \alpha, \beta, \mu, \iota, \Gamma _0) $ or $ ~(\Gamma, \Gamma_{0}).$
 Whenever we write a product in a given a groupoid $~(\Gamma, \Gamma_{0}),~$  we are assuming that it is defined.

A group $~G~$ having $~e~$ as unit element is just a groupoid over $~\{ e\},~$ and conversely, every groupoid with one unit element is a group.

A $~\Gamma_{0}~$-groupoid $~\Gamma~$ is said to be {\it transitive}, if the map\\
 $~\alpha\times\beta: \Gamma\longrightarrow \Gamma_{0}\times \Gamma_{0} , ~x~\to~ (\alpha\times \beta)(x) = (\alpha(x),\beta(x))~$ is surjective.
The map $~\alpha\times \beta~$ is called the {\it anchor} of $~\Gamma.~$\\[-0.5cm]

\markboth{Gheorghe IVAN}{On three generalizations of the group concept .....}

\begin{Ex}
{\rm Let $~{\bf R}^{\ast}~$ be the set of nonzero real numbers and $~a, b\in {\bf R}^{\ast}~$ such that $~a b=1. $ Consider the sets  $ \Gamma = {\bf R}^{\ast}\times {\bf R}^{\ast}: = {\bf R}^{\ast 2}~ $ and $ \Gamma_{0} =\{ (a, a x)\in \Gamma~|~(\forall) x\in  {\bf R}^{\ast}\}. $  It is easy to see that $~ y=ax ~$ if and only if $ x =b y $ for all $~(x, y)\in \Gamma. $ Then
$ (\Gamma, \alpha, \beta,  \odot, \iota, \Gamma_{0})~$ is a groupoid, where the set $~\Gamma_{(2)}~$ and  its structure functions are given by:\\
$\alpha : \Gamma~\rightarrow~ \Gamma_{0}, (x,y)~\rightarrowtail \alpha (x,y):=(x,ax),~\beta : \Gamma~\rightarrow~ \Gamma_{0}, (x,y)~\rightarrowtail \beta (x,y):=(by,y),~$;\\
$ \Gamma_{(2)} = \{ ((x,y), (z, u))\in \Gamma \times \Gamma ~|~ \beta (x,y) =\alpha (z, u)\} = \{ ((x,y), (z, u)) \in \Gamma \times \Gamma ~|~ z=by \};$\\
$(x,y)~\odot~(by,u):=(x, u),~ (\forall)~x, y, u \in {\bf R}^{\ast}~~$ and $~~\iota(x,y):=(by, ax), (\forall)~x, y\in {\bf R}^{\ast}. $
It is easy to verify that the conditions $ (1) - (3) $ of Definition 2.1 hold. First, we consider $~ (x, y), (z, u), (v, w) \in \Gamma.~$   The product
$ (x, y)\cdot (z, u) \cdot (v, w) $ is defined if and only if $~z=by~$ and $~v=bu.~$ We have:\\
$ (1)~~~((x, y)\odot (by, u))\oplus (bu, w)=((x, u)\odot (bu, w) = (x,w)= (x, y)\odot ((by, u)\oplus (bu, w)); $\\
$ (2)~~~\alpha (x, y)\odot (x, y))= (x, ax)\odot (x, y)= (x, t)\odot (bt, y)=(x, y) ~$ and $~(x, y)\odot \beta(x,y)=(x,y)\odot (by, y)= (x,y);$\\
$ (3)~~~(x, y)\odot \iota(x, y)= (x, y)\odot (by, ax)= (x, ax)=\alpha (x,y)~$ and $~\iota(x,y)\odot (x,y)= (by, ax)\odot (x,y)==(by, t)\odot (bt,y)=(by,y)=\beta(x, y).$

This groupoid is denoted with $~{\bf R}^{\ast 2}(a,b).~$ Therefore, for each pair $~(a,b)\in {\bf R}^{\ast 2} ~$   satisfying relation $~ ab=1,~$ an example of groupoid is obtained in this way. We thus have an infinity of groupoids associated with $~{\bf R}^{\ast 2}.~$}\\[-0.6cm]
\end{Ex}
In the following proposition we summarize some properties of the structure functions of a groupoid obtained directly from definitions.
\begin{Prop} {\rm \cite{giv1b})}
~For any groupoid $~(\Gamma, \alpha, \beta, \mu, \iota ,\Gamma_{0})~$ the following assertions hold:\\
$(i)~~~\alpha(u) = \beta(u) =\iota(u)= u ~$ and $~u\cdot u = u~$ for all $~u\in \Gamma_{0};~$\\
$(ii)~~~\alpha(x y) = \alpha(x)~$ and $~ \beta(x y ) = \beta(y), ~(\forall)~ (x,y)\in \Gamma_{(2)};~$\\
$(iii)~~~ \alpha \circ i = \beta, ~~ \beta\circ i = \alpha ~$ and $~ i \circ i = Id_{\Gamma};~$\\
$(iv)~~~$({\it cancellation law} ) If $~x z_{1} = x z_{2}~$ (resp., $~z_{1} x = z_{2} x $), then $~z_{1} = z_{2};~$\\
$(v)~~~(x,y)\in \Gamma_{(2)}~~\Longrightarrow~~(y^{-1},x^{-1})\in \Gamma_{(2)}~~$ and $~~(x y )^{-1} = y^{-1} x^{-1};~$\\
$(vi)~~~$ For each $~u\in \Gamma_{0},~$ the set  $\Gamma(u)=\alpha^{-1}(u)\cap\beta^{-1}(u)=\{ x\in\Gamma~|~\alpha(x)=\beta(x)=u~\}$ is a group under the restriction
of $\mu~$ to $~\Gamma(u),$ called the {\bf isotropy group  of $~\Gamma $ at} $u;$\\
$(vii)~~~$  For each $~x\in \Gamma, $ the isotropy groups $~\Gamma(\alpha(x))~$ and $~\Gamma(\beta(x))~$ are isomorphic;\\
$(viii)~~$ If $~(\Gamma, \Gamma_{0})~$  is transitive, then the isotropy groups $~\Gamma(u),~u\in\Gamma_{0}$ are isomorphes.
\end{Prop}
\begin{defn}
 {\rm  Let $(K,K_{0})$ be a pair of  non-empty subsets, where  $~K\subset \Gamma$ and $ K_{0}\subset \Gamma_{0}.~$\\
 $(i)~$  The pair $(K,K_{0})$  is called {\it  subgroupoid} of  groupoid $~(\Gamma, \Gamma_{0}),~$ if it is closed under multiplication (when it is defined)  and inversion, i.e.
the following conditions hold:\\
$(1)~~~(\forall) ~x,y\in K~$ such that $~x y~$ is defined, we have $~x y\in K;~$\\
$(2)~~(\forall) ~x\in K ~~\Rightarrow ~~x^{-1}\in K.$\\
$(ii)~$ If $~K_{0}=\Gamma_{0},~$ then $~(K,\Gamma_{0})~$ is called  {\it wide subgroupoid} of $~(\Gamma, \Gamma_{0}).~$\\
$(iii)~$ By a {\it normal subgroupoid} of a groupoid $~(\Gamma,  \Gamma_{0}),~$ we mean a wide subgroupoid $~(K, \Gamma_{0})~$ of $~(\Gamma, \Gamma_{0})~$ satisfying the property :\\
 $~$ for all $~x\in \Gamma~$ and $~a\in K~$ such that the product $~x a x^{-1}~$ is defined, we have $~x a x^{-1}\in K.~$}
\end{defn}
If $~\Gamma~$ is a $~\Gamma_{0}~$- groupoid, then $~Is(\Gamma) = \{ x\in \Gamma ~|~\alpha(x)=\beta(x)~\}= \cup_{u\in \Gamma_{0}}\Gamma(u)~$ is a groupoid, called  the
{\it isotropy group bundle }  of  $~\Gamma.~$     It is easy to prove that $~Is(\Gamma) $ is a normal subgroupoid of $~\Gamma,~$ called the  {\it isotropy subgroupoid} of $~\Gamma.~$ 
\begin{defn}
~ {\rm $~(i)~$ Let $~( \Gamma, \alpha, \beta, \mu, \iota , \Gamma_{0})~$ and $~( \Gamma^{\prime}, \alpha^{\prime}, \beta^{\prime}, \mu^{\prime}, \iota^{\prime}, \Gamma_{0}^{\prime})~$  be two groupoids.
A  {\it morphism between these groupoids} is a pair $~(f, f_{0})~$ of maps, where $~f: \Gamma \longrightarrow \Gamma^{\prime}~$ and $~f_{0} : \Gamma_{0}~\to~\Gamma_{0}^{\prime},~$ such that the
following two conditions are satisfied:\\
$(1)~~~ f(\mu(x,y)) = \mu^{\prime}(f(x), f(y)),~~~(\forall) ~(x,y)\in \Gamma_{(2)};~$\\
$(2)~~~ \alpha^{\prime}\circ f = f_{0}\circ \alpha ~$ and $~ \beta^{\prime}\circ f = f_{0}\circ \beta.~$\\
$(ii)~$ A morphism of groupoids $~(f, f_{0})~$ is said to be {\it isomorphism of groupoids}, if $~f~$ and $~f_{0}~$ are bijective maps.}
\end{defn}
\begin{Ex}
{\rm ~{\bf  The disjoint union of a family of groupoids}. If $~\{~\Gamma_{i}~|~i\in I~\}~$ is a disjoint family of groupoids, let $~\Gamma~=~\cup_{i\in I} \Gamma_{i}~$  and $~\Gamma_{(2)} = \cup_{i\in I}\Gamma_{i,(2)}.~$  Here, two elements  $~x,y\in \Gamma~$ may be composed if and only if 
 they lie in the same groupoid $~\Gamma_{i}.~$ This groupoid is called the {\it disjoint union of  groupoids} $~\Gamma_{i}, i\in I,~$ and it is denoted by $~\coprod_{i\in I} \Gamma_{i}.~$ The unit set of this groupoid is $~\Gamma_{0} = \cup_{i\in I} \Gamma_{i,0},~$ where $~\Gamma_{i,0}~$ is the unit set of $~\Gamma_{i}.~$\\[-0.5cm]

In particular, the disjoint union of the groups $~G_{i}, i\in I,~$ is a groupoid, i.e. $~G = \coprod_{i\in I}{G}_{i},~$ which be called the {\it groupoid associated to family of groups} $~{G}_{i}, i\in I.~$ \\
For example, let $~GL(n;{\bf R})~$ the general linear group of order $~n~$ over $~{\bf R}.~$  We consider the family of groups $~\{~ GL(k;{\bf R})~|~ 1\leq k\leq n~\}.~$
Then the disjoint union $~\bigcup\limits_{k=1}^{n} GL(k;{\bf R})~$ is a groupoid, called the {\it general linear groupoid of order }$~n~$ {\it over} $~{\bf R}~$ and denoted by $~{\cal GL}(n;{\bf R}).~$\\
Let the map $~f : {\cal GL}(n;{\bf R})~\to~ {\bf R}^{*}~$ defined by $~f(A) = det (A),~$ for all matrix $~A\in {\cal GL}(n;{\bf R}).~$  We have that $~f~$ is a surjective groupoid morphism from the groupoid $~{\cal GL}(n;{\bf R})~$ onto the multiplicative group of reel numbers different from zero.} 
\end{Ex}
\begin{Ex}
~ {\rm $(i)~$ Let $~M~$  be a nonempty set.  By {\it a quasipermutation of the set} $~M~$ we mean an injective map from a subset of $~M~$ into $~M.~$\\
We denote by $~\Gamma = {\cal S}(M)~$  the set of all quasipermutations of $~M,~$ i.e. $~{\cal S}(M) = \{ f~|~ f:A\to M, f~ \hbox{is injective and } \emptyset \neq A\subseteq M~\}.~$
For $~f\in {\cal S}(M),~$ let D(f) be the domain of $~f, ~R(f)=f(D(f))~$ and\\
$~\Gamma_{(2)} = \{~(f,g)~|~ R(f)=D(g)~\}.~$  For $~(f,g)\in \Gamma_{(2)}~$ we define $~\mu(f,g) = g\circ f.~$\\
If $~Id_{A}~$ denotes the identity on $~A,~$ then $~\Gamma_{0} = \{~ Id_{A}~|~ A\subseteq M~\}~$ is the set of units of $~\Gamma,~$ denoted by $~{\cal S}_{0}(M)~$  and $~f^{-1}~$ is the inverse function from  $~R(f)~$ to $~D(f).~$ The maps $~\alpha,~\beta~$ are defined by  $~\alpha(f) = Id_{D(f)},~ \beta(f) = Id_{R(f)}.~$ Thus $~{\cal S}(M)~$ is a groupoid over $~{\cal S}_{0}(M).~~{\cal S}(M)~$ is called the {\it symmetric groupoid of} $~M~$ or the {\it groupoid of quasipermutations of} $~M.~$\\
$(ii)~$ For a given groupoid $~(\Gamma, \Gamma_{0}),~$ let $~({\cal S}(\Gamma),{\cal S}_{0}(\Gamma))~$ be the symmetric groupoid of the set $~\Gamma,~$ where $~{\cal S}_{0}(\Gamma)=\{ Id_{A}~|~ A\subseteq \Gamma \}.~$\\
We consider now the set $~{\cal L}(\Gamma) = \{~L_{a}~|~a \in \Gamma~\}~$  of all left translations\\
 $L_{a}: \Gamma ~\longrightarrow~\Gamma, ~x \to~ L_{a}(x) = a x,~$ whenever $~(a,x)\in \Gamma_{(2)}.~$\\
We have $~D(L_{a})= \{~x\in \Gamma~|~(a,x) \in \Gamma_{(2)}~\}\neq \emptyset,~$ since $~(a,\beta (a))\in \Gamma_{(2)}~$ and so $~L_{a}\in {\cal S}(\Gamma).~$ Hence  ${\cal L}(\Gamma)~$ is a subset of $~{\cal S}(\Gamma).~$\\
For all $~a,b,x\in \Gamma~$ such that $~\beta(a)=\alpha(b)~$ and $~\beta(b)=\alpha(x)~$ we have\\
 $~L_{a}( L_{b}(x)) = L_{a}(b x) = a(b x) = (a b) x = L_{a b}(x)~$ and we note that $~ L_{a} \circ L_{b} = L_{a b}~$ if $~(a,b) \in \Gamma_{(2)}.~$ Consequently, we have $~ L_{\alpha(x)} \circ L_{x} = L_{x} \circ L_{\beta(x)} = L_{x},~(\forall)~x\in \Gamma.~$\\
For all $~u\in \Gamma_{0}~$ we have $~L_{u} = Id_{D(L_{u})},~$ hence $~L_{u}\in {\cal S}_{0}(\Gamma)~$  and\\
${\cal L}_{0}(\Gamma) = \{ L_{u} ~|~ u\in \Gamma_{0}\}~$ is a subset of $~{\cal S}_{0}(\Gamma).~$ Since $~{\cal L}(\Gamma) \subseteq {\cal S}(\Gamma)~$ and the conditions of Definition 2.2 are satisfied, it follows that $~{\cal L}(\Gamma)~$ is a subgroupoid of $~{\cal S}(\Gamma).~$\\
This groupoid is called the {\it groupoid of left translations of} $~\Gamma.~$}
\end{Ex}
The groupoid of left translations of a given groupoid plays an important role in the proof of Cayley's theorem for groupoids (Theorem 2.1).
\begin{Th} (\cite{giv1b})
Every groupoid $~(\Gamma, \Gamma_{0})~$ is isomorphic to a subgroupoid of the symmetric  groupoid $~({\cal S}(\Gamma), {\cal S}_{0}(\Gamma)).~$
\end{Th}
{\it Proof.} Let $~(\Gamma, \Gamma_{0}) $ be a groupoid and  $~({\cal L}(\Gamma) ,{\cal L}_{0}(\Gamma))~$  its groupoid of left translations. It is proved that the pair of functions $~(\varphi, \varphi_{0}), $
where $~\varphi: \Gamma~ \longrightarrow~ {\cal L}(\Gamma),~\varphi(a) = L_{a},~(\forall)~a\in \Gamma~$ and $~\varphi_{0}: \Gamma_{0} \longrightarrow~ {\cal L}_{0}(\Gamma),~\varphi(u) = L_{u},~(\forall)~u\in \Gamma_{0}~$,
 is an isomorphism of groupoids. \hfill$\Box$\\
\begin{Rem}
{\rm For the detailed presentation of the groupoid concept in the sense of Definition 2.1 and the fundamental properties related to these groupoids,
 as well as their use in different areas of mathematics, we recommend the most relevant references on groupoids (\cite{brown06, dugiv94, higgi71, giv1b, wein96}).
 The basic concepts from groupoid theory  has been defined and studied in a series of papers, for instance:
substructures, normal subgroupoids, quotient groupoid, morphisms of Brandt  groupoids, construction of some ways of building up new groupoids from old ones  (\cite{avmp20, ivan24, mapi22})
 and description of the mportant results on finite groupoids  (\cite{beglpt21, ibr19, miva02, miv3d}).}\\[-0.5cm]
\end{Rem}

                          $~~~~~~~~~~~~~~~~~~~~~~~~~{\bf \star~~~~~\star~~~~~\star~~~~~\star~~~~~\star~~~~~\star~~~~~\star}$ \\ [-0.4cm]                                        

The Brandt groupoids were generalized by C. Ehresmann in \cite{ehre50}. C. Ehresmann added further structures (topological and differentiable as well as algebraic)
to groupoids, thereby introducing them as a tool in topology and differential geometry. After the introduction of topological and differentiable  groupoids by Ehresmann  in the 1950's, they
have been studied by many mathematicians with different approaches (\cite{ehre50, giv1a, mack87, west71}).

In the following we list some fundamental properties  of the concept of groupoid in the sense of Ehresmann which generalize well-known results in group theory (\cite{giv1a, ivan24, miva02}):

${\bf ~-~}$  {\bf definition of the Ehresmann groupoid and presentation of relevant properties} of notions related to it} (types of subgroupoids, normal subgroupoid, groupoid morphism);

 ${\bf ~-~}$  {\bf general construction of Ehresmann groupoids} (direct product of two Ehresmann groupoids;  Whitney sum of two Ehresmann groupoids over the same base; induced groupoid of a Ehresmann groupoid);

${\bf ~-~}$  {i\bf s introduced the concept} of  {\it strong  groupoid morphism}.  This   special type of  groupoid morphism plays an important role in groupoids  theory. More precisely, using the concept of strong groupoid morphism, two theorems are proved, namely: the correspondence theorem for subgroupoids and the  correspondence theorem for normal subgroupoids (\cite{giv1a}). These theorems generalize the correspondence theorems for subgroups and normal subgroups by a surjective morphism of groups.\\[-0.2cm]

{\bf 2.2. Groupoids endowed with topological and differentiable structures}

The research of some geometric objects associated with a topological or differentiable manifold sometimes leads us to a groupoid endowed with a compatible topological or differentiable structure, arriving at the notion of topological groupoid or Lie groupoid.

The concepts of topological groupoid and Lie groupoid are defined using the definition of the groupoid in the sense of Ehresmann.

For basic results  on topological groupoids and bundles of topological groupoids, readers should consult (\cite{brown06, ehre50, ivan24, miva05}).

The study of the basic topics related to Lie groupoids and bundles of Lie groupoids has been addressed in a series of works, among which we mention (\cite{ivan24, miv2a,  miva06,   mack87}).\\[-0.2cm]

{\bf 2.3  Structured groupoids}

Another approach to the notion of a groupoid is that of a structured groupoid.  This  concept is obtained by adding another algebraic structure
such that the composition of the groupoid is compatible with the operation of the added algebraic structure. The most important types of structured groupoids
 are  the following: group-groupoid  (\cite{brspen, ivan24}), vector groupoid (\cite{popu4a}), topological group-groupoid and Lie group-groupoid (\cite{fana13, ivan24}).

The notion of group-groupoid was defined by R. Brown and Spencer in \cite{brspen}.  An extension of the algebraic concept of group-groupoid to the notion of vector-space groupoid, was defined and investigated by Mihai Ivan in the paper 
\cite{miva13}.

  The group-groupoids, vector-space groupoids and their extensions(topological group-groupoids,  Lie group-groupoids) are mathematical structures that have proved to be useful in many areas of science, see for instance 
the papers (\cite{akiz18, guari17, guric18}).\\[-0.2cm]

{\bf  2.4. Applications of  groupoids}

A large number of researchers have paid great attention to the study of groupoids and their applications. These researchers have investigated various topics related to groupoids, achieving remarkable  results in  mathematics (algebra, analysis, geometry)
and science,  see (\cite{avmp20, beglpt21, ivan24, rare01, wein96, west71}).

The study of Ehresmann groupoids, topological groupoids, Lie groupoids and Lie algebroids is motivated by their applications in many  branches of mathematics and engineering, namely in:\\
- category theory, topology, representation theory, algebraic topology (\cite{brspen, brown06, ehre50, higgi71, ibr19, igim11});\\
- differential geometry, nonlinear dynamics, quantum  mechanics (\cite{goste06, puta5a, gmio06, opris5a,igmp11, mack87});\\
- study of some dynamical systems in Hamiltonian form and fractional dynamical systems using Lie groupoids and Lie algebroids (\cite{giom13, ivmi16, migi18, mlag, miva23, miva24, migo09}).

Moreover, various types of constructions and classes of actions of Lie groupoids on differential manifolds are used in the construction of principal bundles with structural Lie groupoids (\cite{giv5b, giv5c}). These principal bundles constitute an important 
study tool in various applications of differential geometry in theoretical mechanics. \\[-0.6cm]

\section{Generalized groups: definitions, generalized subgroup and generalized homomorphism}

In this section we give the defnition of generalized group and some of its basic properties. The generalized group is an interesting generalization of groups. The concept of a generalized group was introduced by R. M. Molaei  ( \cite{mola98}). 
While there is only one identity element in a group, each element in a generalized group has a unique identity element. 
\begin{defn} (\cite{mola99})~ {\rm  A {\it generalized group} is a nonempty set $~{\cal G}~$ together with a binary operation (called {\it multiplication})  $~ m:\ {\cal G}\times {\cal G}\rightarrow {\cal G},~(a,b)\longmapsto  m ( a,b):=a b,~$  satisfying the following conditions:\\ 
$(1)~$({\it associativity})~$ (a b) c = a (b c),~$  for all  $~a, b, c \in {\cal G}$ ;\\
$(2)~$({\it existence and uniqueness of identity}): for each $~a\in {\cal G}~$ there exists a unique element $~e(a)\in {\cal G}~$  such that $~a e(a) = e(a) a = a $;\\
$(3)~$({\it existence of inverse element}): for each $~a\in {\cal G}~$ there exists  $~a^{-1} \in  {\cal G}~$  such that $~a a^{-1} = a^{-1} a =e(a). $}
\end{defn}

In the following proposition  we give some characteristic properties related to the structure of generalized group.

\begin{Prop} {\rm \cite{mola98})} Let  $~{\cal G}~$  be a generalized group. Then  the following assertions hold:\\
$(i)~~~$ for each  $~a \in {\cal G}~$  there exists a unique element $~a^{-1}\in  {\cal G};$\\
$(ii)~~~$ for each $~a \in {\cal G}, $ we have $~ e(a) = e(a^{-1})~$  and $~e(e(a)) = e(a);$\\
$(iii)~~~$  for each $~a \in {\cal G}, $ we have $~(a^{-1})^{-1}=a .$
\end{Prop}

We notice that  for each element $~a~$  in a generalized group  $~{\cal G},$ the inverse  $~a^{-1}~$ is unique and both $~a~$ and  $~a^{-1}~$ have the same identity.

 Let  $~{\cal G}~$  be a generalized group.  It is called {\it abelian} or {\it commutative generalized group~} if $~ a b= b a~$ for all  $~a, b \in {\cal G}.$

\begin{Rem} {\rm  Let  $~G~$  be a  group having  $~e_{0}~$ as identity element. Then   $~G~$ is a generalized group, taking $~e(a):=e_{0}~$ for all $~a\in G$.}
\end{Rem}

By Remark 3.1, we  can conclude  that every group is a generalized group, but it is not true in general that every generalized group is a group.

Proposition 3.2 states the relationship between groups and generalized groups.

\begin{Prop} {\rm \cite{guari17, mola99})}Let  $~{\cal G}~$  be a commutative generalized group. Then $~{\cal G}~$ is a group.
\end{Prop}
\begin{defn} (\cite{mola99})~ {\rm  A generalized group $~{\cal G}~$  is called {\it normal generalized group} if $~ e(ab) = e(a)e(b)~$ for all  $~a, b \in {\cal G}.$}
\end{defn}
\begin{defn} (\cite{mola99})~ {\rm  A  non-empty subset $~{\cal H}~$ of a generalized group  $~{\cal G}~$ is called  a {\it generalized subgroup} of  $~{\cal G},~$ if  $~a b^{-1} \in {\cal H}~$ for all $~a, b\in {\cal H}.$ }
\end{defn}
\begin{Prop} {\rm Let $~({\cal H}_{i})_{\in I}~$ be a family of generalized subgroups of the generalized group  $~{\cal G}~$  and $~\cap {\cal H}_{i}\neq \emptyset.~$  Then  $~\cap {\cal H}_{i}, i\in I~$  is a generalized subgroup of $~{\cal G}.$}
\end{Prop}
\begin{Prop} {\rm Let $~{\cal G}~$ be a  generalized group and $~a\in {\cal G}.~$  Then  $~{\cal G}_{a}:=\{   x \in  {\cal G}~|~e(x)=e(a) \} $ is a generalized subgroup of $~{\cal G}.$}\\[-0.5cm]
\end{Prop}
{\it Proof.} The proof of this proposition can be found in paper \cite{mem00}.\hfill$\Box$
\begin{Rem} {\rm  Let $~{\cal G}~$ be a  generalized group. For each  $~a\in {\cal G}, ~$  it  follows that $~{\cal G}_{a}~$ is a group having $~e(a)~$ as identity element}
\end{Rem}

\begin{Ex} ~ {\rm  Let ${\cal G} =\{ A=\left (\begin{array}{cc}

a_{1}  & a_{2}\\

a_{3} & a_{4}

\end{array}\right)~|~a_{i} \in {\bf R}, i=\overline{1,4}~ \hbox{and}~ det(A)>0 \}.~ $ It easy to prove that for all $~\alpha\in {\bf R}, \alpha>0~$ and  $~A\in {\cal G},~$ we have:
 $\sqrt{det(\alpha\cdot A)} = \alpha\cdot \sqrt{det(A)}.$ \\
 We now define a generalized  group  structure on the set $~{\cal G}.~$ Define on $~{\cal G}~$ a multiplication $~\odot: {\cal G}\times {\cal G}~ \rightarrow ~{\cal G}~$ by $~(A,B) \rightarrow A\odot B,~$ where:
$~A\odot B:= \sqrt{det(A)} \cdot B.$

For all $~X, Y \in {\cal G},~$ the following relation  holds:\\[0.2cm]
$~~~~~~~~~~~~~~~~~~~~~ \sqrt{det( X\odot Y)} =\sqrt{det(X)}\cdot \sqrt{det(Y)}.$

For all $~A\in {\cal G},~$ the identity element  of  $~A~$ and the inverse of $~A~$ are given by:\\[0.2cm]
$~~~~~~~~~~~~~~~~~~~~~e(A) = \frac{1}{\sqrt{det(A)}}  \cdot A~~~~~$  and  $~~~~~A^{-1} = \frac{1}{det(A)}\cdot A.$\\[0.2cm]

It is easy to verify that the conditions $ 1-3 $ of Definition 3.1 hold. For this, we consider $~A, B, C\in {\cal G}~$. We have:\\[0.2cm]
$ (1)~~~(A\odot B)\odot C=(\sqrt{det(A\odot B})\cdot C=(\sqrt{det(A)}\cdot \sqrt{det(B)})\cdot C~~$ and $~~A\odot (B\odot C)=\sqrt{det(A)}\cdot (\sqrt{det(B\odot C)}=\sqrt{det(A)}\cdot (\sqrt{det(B)}\cdot C).~$
Hence, $~(A\odot B)\odot C=A\odot (B\odot C).$\\[0.2cm]
$ (2)~~~e(A)\odot A=\sqrt{det (e(A))}\cdot A= \sqrt{det ( \frac{1}{\sqrt{det(A)}}  \cdot A)}\cdot A = \frac{1}{\sqrt{det(A)}}\cdot \sqrt{det(A)}\cdot A=A~~$ and \\[0.2cm]
$~A\odot e(A)=\sqrt{det(A)}\cdot e(A)= \sqrt{det(A)}\cdot \frac{1}{\sqrt{det(A)}}  \cdot A =A.~$  Hence  $~e(A)\odot A = A\odot e(A)=A.$
$(3)~~~A\odot A^{-1}=\sqrt{det (A)}\cdot A^{-1}=\sqrt{det (A)}\cdot  \frac{1}{det(A)}\cdot A = \frac{1}{\sqrt{det(A)}}\cdot A= e(A)~~~$  and\\[0.2cm]
$~A^{-1}\odot A=\sqrt{det (A^{-1})}\cdot A=\sqrt{det(\frac{1}{det(A)}\cdot A)}\cdot A=\frac{1}{det(A)}\cdot \sqrt{det(A)}\cdot A=  \frac{1}{\sqrt{det(A)}}  \cdot A=e(A).~$ Hence $~A\odot A^{-1}=A^{-1}\odot A= e(A).$\\[0.2cm]
We conclude that $~{\cal G}~$ is a generalized group.

$\bullet~~~$   {\bf Calculation exercise in the generalized group} $~{\cal G}.$ Consider the matrices  $~A_{1}=\left (\begin{array}{cc}

-1  & 1\\

-2 & -2

\end{array}\right)~ $ and  $~A_{2}=\left (\begin{array}{cc}

1  & -2\\

3 & 3

\end{array}\right). ~$  We have:\\
$~det(A_{1})=4,~det(A_{2})=9,~~~A_{1}\odot A_{2} =2 A_{2} ==\left (\begin{array}{cc}

2  & -4\\

6 & 6 
\end{array}\right).~ $   Then $~e(A_{1})=\frac{1}{2}\cdot A_{1},$\\
$~e(A_{2})=\frac{1}{3}\cdot A_{2}, ~ A_{1}^{-1}=\frac{1}{4}\cdot A_{1} ~$  and  $~ A_{2}^{-1}=\frac{1}{9}\cdot  A_{2}~.$ Also, 
$~A_{2}\odot A_{1} =3 A_{1} =\left (\begin{array}{cc}

-3  & 3\\

-6 & -6 
\end{array}\right)~ $ and we observe that $~{\cal M}~$ is not a commutative generalized group. }
\end{Ex}

\begin{defn} (\cite{mola98})~ {\rm  Let  $~{\cal G}~$ and  $~{\cal G}^{\prime}~$ be two generalized groups. A {\it generalized group homomorphism} from    $~{\cal G}~$ to  $~{\cal G}^{\prime}~$ is a map  $~f: {\cal G}~\rightarrow~{\cal G}^{\prime}~$
such that $~ f(a b) = f(a) f(b)~$ for all $~a,b\in {\cal G}.$}
\end{defn}
\begin{Prop} {\rm Let  $~f: {\cal G}~\rightarrow~{\cal G}^{\prime}~$ be a  generalized group homomorphism.  Then:\\
$(i)~~~f(e(a)) = e(f(a))~$  is an identity element in $~{\cal G}^{\prime}~$ for all $~a\in {\cal G}.$ \\ 
$(ii)~~~f(a^{-1}) = (f(a))^{-1}~$ for all   $~a\in {\cal G}.$\\
$(iii)~~~$  if $~{\cal H}~$   is a generalized subgroup of $~{\cal G},~$  then $~f({\cal H})~$ is a generalized subgroup of $~ {\cal G}^{\prime}.$\\
$(iv)~~~$  if $~{\cal H}^{\prime}~$   is a generalized subgroup of $~{\cal G}^{\prime}~$ and $~f^{-1}({\cal H}^{\prime})\neq \emptyset,~$ then $~f^{-1}({\cal H}^{\prime})~$ is a generalized subgroup of $~ {\cal G}.$} \\[-0.5cm]
\end{Prop}
{\it Proof.} The proof of this proposition can be found in paper \cite{mem00}.\hfill$\Box$

\begin{Ex} ~ {\rm  Let ${\cal M} =\{ A=\left (\begin{array}{ccc}

a  & b & c\\

0 &  1 & 0\\

0 &  0 & 0\\

\end{array}\right)~|~a, b, c\in {\bf R} \}.~ $  Let $~A,~B~$ be two matrices in $~{\cal M}.$  We have:\\

$ A\cdot B= \left (\begin{array}{ccc}

a  & b & c\\

0 &  1 & 0\\

0 &  0 & 0\\

\end{array}\right)  \left (\begin{array}{ccc}

a^{\prime}  & b^{\prime} & c^{\prime}\\

0 &  1 & 0\\

0 &  0 & 0\\

\end{array}\right) =  \left (\begin{array}{ccc}

a a^{\prime}  &a  b^{\prime}+b & a c^{\prime}\\

0 &  1 & 0\\

0 &  0 & 0\\

\end{array}\right) \in {\cal M}.$\\

Therefore, the usual matrix multiplication defines a binary operation on the set $~{\cal M}.$

For all $~ A=\left (\begin{array}{ccc}

a  & b & c\\

0 &  1 & 0\\

0 &  0 & 0\\

\end{array}\right) \in {\cal M},~$ the identity element  of  $~A~$ and the inverse of $~A~$ are given by:\\

$~~~~~~~~~~~~~~~~~~~~~e(A) =\left (\begin{array}{ccc}

1  & 0 & c/a\\

0 &  1 & 0\\

0 &  0 & 0\\

\end{array}\right) ~~~~~  and  ~~~~~A^{-1} =\left (\begin{array}{ccc}

1/a  & -b/a & c/a^{2}\\

0 &  1 & 0\\

0 &  0 & 0\\

\end{array}\right) .$

It is easy to verify that the conditions $ 1-3 $ of Definition 3.1 hold.  We conclude that $~{\cal M}~$ is a generalized group.\\ 

$\bullet~$  Consider $~A, B\in {\cal M}.~$  By direct calculation it is easy to verify that  $~e(A\cdot B)= e(A)\cdot e(B)~$ and it is deduced that $~{\cal M}~$ is a normal generalized group.\\

$\bullet~$  Let be the matrices $~C =\left (\begin{array}{ccc}

1  & -1 & 2\\

0 &  1 & 0\\

0 &  0 & 0\\

\end{array}\right) ,~  D=\left (\begin{array}{ccc}

2  & 1 & -1\\

0 &  1 & 0\\

0 &  0 & 0\\

\end{array}\right) .~$  We have:\\[0.2cm] $~C\cdot D =\left (\begin{array}{ccc}

2  & 0 & -1\\

0 &  1 & 0\\

0 &  0 & 0\\

\end{array}\right) ~$  and  $ ~D\cdot C=\left (\begin{array}{ccc}

4  & 2 & -2\\

0 &  1 & 0\\

0 &  0 & 0\\

\end{array}\right) ~.$ Then $~C\cdot D \neq D\cdot C.~$  Therefore  $~{\cal M}~$  is not a commutative generalized group.

$\bullet~$  Consider the multiplicative group $~{\bf R}^{\star}={\bf R}\setminus \{0\}. $  It is a generalized group. 
Define the map    $~f: {\cal M}~\rightarrow~{\bf R}^{\star}~$ by $~f(A):=a~$ for all $~A\in {\cal M}.$ This map is a generalized group homomorphism. Indeed, for all $~A,B\in {\cal M}~$  we have $~f(A\cdot B) = a\cdot a^{\prime},~f(A)\cdot f(B) = a\cdot a^{\prime}~$ and 
 $~f(A\cdot B)= f(A)\cdot f(B).$}
\end{Ex}
\begin{Rem} {\rm  Generalized groups have been studied by many mathematicians (\cite{guari17, mem00, mola98}). This mathematical concept has been applied in geometry (\cite{mofar09}) and dynamical systems (\cite{mola05}). For the detailed presentation of the generalized group concept and the fundamental properties related  to this generalization of the notion of group, one can consult the works (\cite{mem00, mola98, mola99, nade19}). This concept has been used in defining and studying mathematical concepts in algebraic, topological and differentiable terms that have significant applications in various branches of mathematics, namely: generalized ring (\cite{mola02}), generalized crossed module (\cite{guari17},) topological generalized group (\cite{guari17}),  generalized Lie group (\cite{fana13}) and the construction of some ways of building up new groupoids from old ones (\cite{nade19}).}\\[-0.5cm]
\end{Rem}

\section{Almost groupoids}
\indent
In this section we present a new generalization of the notion of group, namely: the {\it almost groupoid} concept. This concept was introduced by Mihai Ivan in paper (\cite{miva23a}).\\[-0.2cm]

{\bf 4.1. Basic properties: almost groupoids and almost subgroupoids}

We will give the definition of an almost groupoid from a purely algebraic point of view.
\begin{defn}(\cite{miva23a})
{\rm   An {\it almost groupoid $ G$ over} $ G_{0}$   is a pair
$( G, G_{0} )$ of nonempty sets such that $~G _0\subseteq G,~ $ endowed with a surjective map
$\theta: G \rightarrow G_{0}$, a partially binary operation $~ m : G_{(2)}\rightarrow G ,~(x,y)\longmapsto  m ( x,y):=x\cdot y,~$ where\\
$~~~~~~~~~~~~~~~~~~~~~~~~~~~~ G_{(2)}:= \{ (x,y)\in G\times G | \theta(x) = \theta(y) \}$ \\
and a map $\iota:G \rightarrow G ,~~x\longmapsto \iota(x):=x^{-1},$ satisfying the following properties:\\
$({\bf AG1})~$ ({\it associativity})$~~(x\cdot y)\cdot z=x\cdot (y\cdot z)~$ in the sense that if
one of two products $ (x\cdot y)\cdot z $ and $ x\cdot (y\cdot z) $ is defined,
then the other product is also defined and they are equals;\\
$({\bf AG2})~$ ({\it units}): for each
$ x\in G~\Rightarrow~(\theta(x),x),~(x,\theta( x))\in G _{(2)}$
and  we have\\ $ \theta(x)\cdot x= x \cdot \theta(x)= x;$\\
$({\bf AG3})~$ ({\it inverses}): for each $ x\in G~\Rightarrow~(x, x^{-1}),~( x^{-1}, x) \in G_{(2)} $
 and we have\\ $ x\cdot x^{-1}= x^{-1}\cdot x= \theta(x).$}
\end{defn}
\markboth{Gh. Ivan}{On three generalizations 9f the concept group.....}
If $ G $ is an almost groupoid over $ G_{0}, $  we will sometimes write $x\cdot y $ or $ x y$ for $ m(x,y). $ Also, the set $ G_{(2)} $ is called  the {\it set of composable pairs} of $ G. $

 An almost groupoid $ G $ over $ G_{0} $  with the \textit{structure functions} $ \theta
$ ({\it units map}),$ m $ ({\it multiplication}) and $ \iota $ ({\it inversion}), is denoted by $ (G, \theta,  m, \iota, G_{0}).~$  Whenever we write a product in a given almost groupoid, we are
assuming that it is defined.

In view of Definition 2.1, if $ x,y, z\in G, $ then:

 $(i)~~$ the product $~~x\cdot y\cdot z := (x\cdot y)\cdot z ~~$ {\it is defined} $~~~ \Leftrightarrow ~~~ \theta(x)=\theta (y) =\theta (z). $

$(ii)~~ \theta(x)\in G_{0} $ is {\it the unit} of $ x $ and  $ x^{-1}\in G $ is the {\it inverse} of $ x.$\\

For an almost groupoid $ G $  over $ G_{0}, $  we will use the following notations:  $ (G, \theta, m, \iota, G_{0} ) $
or $ (G, G_{0} ). $  Furthermore,  $ G_{0} $ is called the {\it units set} of $ G. $

The basic properties of the almost groupoids are given in the following two propositions.
\begin{Prop}   (\cite{miva23a}) If $ (G,\theta,  m, \iota, G_{0}) $ is an almost groupoid, then:

$(i)~~~\theta(u)=u,~(\forall)u\in G_{0}.$

$(ii)~~u\cdot u=u ~$ and $~\iota (u)=u,~(\forall)u\in G_{0}.$

$(iii)~\theta(x\cdot y) =\theta (x), ~(\forall ) ( x,y ) \in G _{(2)}.$

$(iv)~\theta ( x^{-1}) =\theta ( x)~$ and $~\theta(\theta(x))= \theta(x), ~~  (\forall) x\in G.$

$(v)~ $  if $~(x,y)\in G_{(2)},$ then $~(y^{-1},x^{-1})\in G_{(2)}~ $ and $~(x\cdot y)^{-1}=y^{-1}\cdot x^{-1}.$

{\it(vi)}$~(x^{-1})^{-1}=x, ~~(\forall) x\in G.$

$(vii)~$  If $~ x,y\in G,~$  then:  $~~ x\cdot y ~~$ is defined $~~~\Leftrightarrow~~~ y\cdot x~~$  is defined.

$(viii)~ (\forall)~ a\in G~~\Rightarrow ~~$ the products $~ a^{2}:= a\cdot a, ~ a^{3}:= a^{2}\cdot a~ $ are defined
\end{Prop}

\begin{Prop} Let $ (G,\theta,  m, \iota, G_{0}) $  be an almost groupoid. For each $ u \in G_{0}, $ the set $ G(u):=\theta^{-1}(u)= \{ x\in G | \theta(x)= u \} $ is
a group.
\end{Prop}
The group $ G(u) $ for $ u\in G_{0} $ is called the {\it isotropy group} of $ G $ at $ u.$
\begin{defn}
{\rm   An {\it almost groupoid $ G$ over} $ G_{0}$ is called {\it abelian} or {\it commutative} if the isotropy group $ G(u) $ is abelian for all $ u\in G_{0}. $}
\end{defn}
\begin{Ex}
{\rm $(i)~$ A group $ G $ having $ e $ as unity, is an almost groupoid over   $ \{e
\}$ with respect to structure functions: $~\theta (x): = e, ,(\forall)x\in G; ~G_{(2)}= G\times G,~ m (x,y):= xy,~ (\forall) x,y\in G
  $ and $ \iota:G \to G,~\iota(x):= x^{-1}, (\forall)x\in G.$

$(ii)~$ A nonempty set $ G_{0} $ may be regarded to be an
almost groupoid over $ G_{0}, $ called the {\it null almost groupoid} associated
to $ G_{0}. $ For this, we take $~ \theta=\iota = Id_{G_{0}} $ and\\
 $ u\cdot u = u, ~(\forall)  u\in G_{0}.$}
\end{Ex}
\begin{Ex} {\rm  Let $ B_{2}=\{ 0, 1\}~$ a set consisting of two elements. Consider the group $~({\bf Z}_{n}, + ) $ of integers modulo $ n, $ where $~{\bf Z}_{n} = \{ c_{1}=0,  c_{2}=1,\ldots, c_{n-1}=n-2,  c_{n}=n-1 \}.~$  We now define an almost groupoid  structure on
 the set $ G:=B_{2}\times  {\bf Z}_{n}~.$ Let $ G_{0} = \{ (a, c_{1}) | a\in B_{2}\}. $  Then $~ (G, \theta, \oplus, \iota, G_{0})~$ is an almost groupoid over $ G_{0}, $ where  the set of composable pairs $~ G_{(2)}~$ and  its structure functions are given by:\\
$\theta : G=B_{2}\times  {\bf Z}_{n}~ \rightarrow~ G_{0}, (a, c_{j})~\rightarrowtail \theta (a, c_{j}):=(a, c_{1});$\\
$ G_{(2)} = \{ ((a, c_{j}), (b, c_{k}))\in G^{2}~|~\theta (a, c_{j}) =\theta (b, c_{k})\} = \{ ((a, c_{j}), (b, c_{k}))\in G^{2}~|~ a=b \};$\\
$(a, c_{j})~\oplus~(a, c_{k}):=(a, c_{j}+ c_{k}),~ (\forall)~a \in B_{2}~, c_{j}, c_{k}\in {\bf Z}_{n}~$ and\\
$\iota(a, c_{j}):=(a, -c_{j}), (\forall)~a\in B_{2},~c_{j}\in  {\bf Z}_{n}~.$  
It is easy to verify that the conditions $ ({\bf AG1}) - {\bf AG3}) $ of Definition 4.1 hold.\\
First, we consider $~ x, y, z\in G ,~$ where $~x=(a_{1}, c_{i}), y=(a_{2}, c_{j}), z=(a_{3}, c_{k}). $  The product $ x\oplus y \oplus z $ is defined if and only if $~a_{1} = a_{2} = a_{3}.$ We have:\\
$ (1)~~~(x\oplus y)\oplus z= ((a_{1}, c_{i})\oplus (a_{1}, c_{j}))\oplus (a_{1}, c_{k})= (a_{1}, c_{i}+c_{j})\oplus (a_{1}, c_{k})= (a_{1}, c_{i}+ c_{j}+ c_{k}) =x\oplus (y\oplus z); $\\
$ (2)~~~x\oplus \theta(x)= (a_{1}, c_{i})\oplus \theta(a_{1}, c_{i}))= (a_{1}, c_{i})\oplus(a_{1}, c_{1})=(a_{1}, c_{i}+0)= x =\theta(x)\oplus x;$\\
$ (3)~~~x\oplus \iota(x)= (a_{1}, c_{i})\oplus \iota(a_{1}, c_{i}))= (a_{1}, c_{i})\oplus(a_{1}, -c_{i})=(a_{1}, 0)= \theta(x)= \iota(x)\oplus x.$\\
We conclude that $ G = B_{2}\times {\bf Z}_{n}~$ is an almost groupoid over $ G_{0}=\{ (0, c_{1}),~(1, c_{1})\}.~$  Observe that $~ (G, G_{0})~$ is a finite almost groupoid of order $~2n.$ }
\end{Ex} 
\begin{defn}
{\rm  Let $(G ,\theta, m, \iota, G_{0})$ be an almost groupoid. A pair of nonempty subsets $(H,H_{0})$
where $ H\subseteq G $ and ${H_0}\subseteq G_{0} $, is called
\textit{almost subgroupoid} of $ G, $ if:

$(1.1)~~\theta (H)=H_{0};$

$(1.2)~~H $ is closed under multiplication and inversion,
that is:\\[0.1cm]
$(1.2.1)~(\forall)~x,y\in H$  such that $(x,y)\in G_{(2)}
\Longrightarrow ~x\cdot y\in H;$\\[0.1cm]
 $~(1.2.2)~(\forall)~x\in
H\Longrightarrow x^{-1}\in H.$}
\end{defn}
\begin{defn}
{\rm Let $(H,H_{0})$ be an almost subgroupoid of the almost groupoid $ (G,G_{0}).$

{\it(i)} $ (H, H_{0}) $ is said to be a {\it wide almost subgroupoid} of $ (G, G_{0} ), $ if $~H_{0}= G_{0},$   that is  $ H $  and $ G $  have the same units.

{\it(ii)} A wide almost subgrupoid $ (H,H_{0})$  of $ (G,G_{0}),$  is called {\it normal almost subgroupoid,} if for all $ g\in G $ and for all
 $ h \in H $ such that the product $ g\cdot h\cdot {g^ {-1}} $ is defined, we have $ g\cdot h\cdot {g^ {-1}}\in H .$}
\end{defn}
\begin{Ex}
{\rm $~(i)~$ Each subgroup (resp., normal) of the group $ ~(G, \{e\}) ~$  is a wide almost subgroupoid (resp., normal) of the almost groupoid $ ~(G, \{e\}). ~$\\
$~(ii)~$ $ G_{0} $ is a normal almost subgroupoid of the almost groupoid  $ (G, G_{0}). $ }
\end{Ex}
\begin{Prop} (\cite{miva23a}) Let $ (G, G_{0}) $  be an almost groupoid.  Then:\\
$(i)~~$ The isotropy group  $~G(u), u\in G_{0}~$ is an almost subgroupoid of $ G. $\\
$(ii)~~$  The set $~Is(G) = \cup_{u\in G_{0}}G(u)\subset G~$ is a normal almost subgroupoid of $~G.~$
\end{Prop}
The almost groupoid $~(Is(G), G_{0}) $ is called the {\it isotropy almost subgroupoid} of $~G. $\\[-0.2cm]

{\bf 4.2. Morphisms of almost groupoids}
\begin{defn}
{\rm  {\it (i)~} Let $ (G, \theta, m, \iota, G_{0}) $ and $ (G^{\prime}, \theta^{\prime}, m^{\prime}, \iota^{\prime}, G_{0}^{\prime}) $ be two
almost groupoids. A {\it morphism of almost groupoids} or {\it almost groupoid
morphism} from $ (G, G_{0}) $ into $ (G^{\prime}, G_{0}^{\prime}) $ is a pair $(f, f_{0}) $
of maps, where $ f:G \to G^{\prime} $ and $ f_{0}: G_{0} \to
G_{0}^{\prime} $ such that the following conditions hold:

$(1)~~~f(m(x,y)) = m^{\prime}(f(x),f(y))~~~\forall (x,y)\in
G_{(2)};$

$(2)~~~\theta^{\prime}\circ f = f_{0} \circ
\theta. $

{\it (ii)~} An almost groupoid morphism $ (f, f_{0}): (G, G_{0}) \to (G^{\prime},
G_{0}^{\prime})$ such that $ f $ and $f_{0}$ are bijective maps, is
called {\it isomorphism of almost groupoids}.

{\it (iii)~}  If  $ (f, f_{0}): (G, G_{0}) \to (G^{\prime}, G_{0}^{\prime})$ is an isomorphism,  
 we say that $ G~$ and $~G^{\prime}~$  are isomorphic almost groupoids and we will write $~ G \cong  G^{\prime}.$}
\end{defn}

\begin{Ex} {\rm   Let $~ G = B_{2}\times {\bf Z}_{n}~$ be  the  almost groupoid given in Example 4.2. Consider the  group $~({\bf Z}_{n}, + ) .$ It is an almost groupoid. 
Define the map    $~f:  B_{2}\times {\bf Z}_{n}~\rightarrow~{\bf Z}_{n}~$ by $~f((a,c_{j})):=c_{j}~$ for all $~a\in B_{2}~$ and $~c_{j}\in {\bf Z}_{n}.$ This map is a morphism of almost groupoids. Indeed, for all $~(a, c_{j}), (a, c_{k}) \in B_{2}\times {\bf  Z}_{n}~$  
we have $~ f((a, c_{j})~\oplus~(a, c_{k})) =f((a, c_{j}+ c_{k}))= c_{j}+ c_{k}~$   and $~  f((a, c_{j}))+ f((a, c_{k})) =  c_{j}+ c_{k}.~$  Hence, $~ f((a, c_{j}))~\oplus~(a, c_{k})) = f((a, c_{j})) + f((a, c_{k})).$ }
\end{Ex}
\begin{Rem}~{\rm For important results regarding the morphisms of  almost groupoids, see \cite{miva23a}.}\\[-0.7cm]
\end{Rem}
\begin{Rem} {\rm ~{\it  The pair groupoid over a set}. Let $~X~$ be a nonempty set. Then $~ \Gamma = X \times X~$ is a groupoid with respect to rules:  $~\alpha(x,y) = (x,x),~ \beta(x,y) = (y,y),~$ the elements $~ (x,y)~$ and $~(y^{\prime},z)~$
are composable in $~\Gamma~\Leftrightarrow~y^{\prime} = y~$ and we take $~(x,y)\cdot (y,z) = (x,z)~$ and the inverse of $~(x,y)~$ is defined by $~(x,y)^{-1} = (y,x).~$\\
This groupoid  is called the {\it pair groupoid} associated to $~X~$ and it is denoted with $~{\cal PG}(X).~$ The unit space of  $~{\cal PG}(X)~$  is the diagonal $~\Delta_{X} = \{ (x,x)~|~ x\in X\},~$  denoted with   $~{\cal PG}_{0}(X).$
  The groupoid $~({\cal PG}(X), {\cal PG}_{0}(X))~$ {\bf is not an almost groupoid}, because $~\alpha(x,y)=(x,x)\neq \beta(x,y)=(y,y) ~$ for  $ x\neq y.$}\\[-0.2cm]
\end{Rem}

{\bf 4.3. Construction of new almost groupoids}

$\bullet~~~$   {\bf  The disjoint union of two  almost groupoids }

Let  $~(G, \theta, \mu, \iota, G_{0}) $  and  $~(G^{\prime}, \theta^{\prime}, \mu^{\prime}, \iota^{\prime}, G_{0}^{\prime}) $ be two almost groupoids such that \\
$~G\cap G^{\prime}= \emptyset.$ 
On the set $~ G\coprod G^{\prime}:= G\cup  G^{\prime}~$  we define  a structure of almost groupoids as follows.

Consider the  set of composable pairs   $~G_{(2)}\cup G_{(2)}^{\prime}.~$  Here, two elements  $~x,y\in G \cup G^{\prime}~$ may be composed if and only if  they lie in the same almost groupoid $~G~$ or $~G^{\prime}.$

This almost groupoid is called the {\it disjoint union of  almost groupoids} $~G~$ and  $~G^{\prime}.~$  It is denoted by $~G\coprod G^{\prime}~$  and its  unit set is $~G_{0}\cup G_{0}^{\prime}.$

$\bullet~~~$   {\bf The direct product  of two  almost groupoids } 

Let  $~(G_{k}, \theta_{k}, \mu_{k}, \iota_{k}, G_{k,0})~$  for $~k= 1, 2~$ be two almost  groupoids. Consider\\
$~G = G_{1}\times  G_{2}~$ the direct product of the sets $~G_{k},~k=1,2.~$ We give to $~G~$ a structure of  almost groupoid as follows.
 The elements $~x = (x_{1}, x_{2})~$ and $~y = (y_{1}, y_{2})~$ from $~G~$ are composable if and only if $~(x_{k},y_{k})\in G_{k,(2)}~$ for $~k=1,2 $ and
we take\\
 $~(x_{1}, x_{2} )\cdot (y_{1}, y_{2})= ( x_{1} y_{1}, x_{2} y_{2} ).~$ \\
It is easy to verify that $~(G_{1}\times G_{2},\theta_{1}\times\theta_{2},  \mu_{1}\times \mu_{2}, \iota_{1}\times \iota_{2}, G_{1,0}\times G_{2,0})~$ 
is an almost  groupoid, called the {\it direct product of almost groupoids} $~(G_{1},  G_{1,0})~$ and  $~(G_{2}, G_{2,0}).~$

$\bullet~~~$   {\bf The semidirect product  of two  almost groupoids } 

\begin{defn}{\rm  Let $ (G, \theta_{1}, \otimes_{1}, \iota_{1}, G_{0}) $ and  $ (H, \theta_{2}, \otimes_{2}, \iota_{2}, H_{0})~$ be two almost groupoids. An {\it almost groupoid action} of the almost groupoid $~G~$ 
on the almost groupoid $~H~$  is a map $~\odot: G \times H \rightarrow H,~ (g,h)~\longmapsto g\odot h~$  such that the following conditions are satisfied:\\
$(1)~ (g_{1}\otimes_{1} g_{2})\odot h = g_{1}\odot (g_{2}\odot h),~~ \forall~(g_{1}, g_{2}) \in G_{(2)} ~$ and $~h\in H.$\\
$(2)~~g\odot (h_{1}\otimes_{2} h_{2})= (g\odot h_{1}))\otimes_{2} (g\odot h_{2}),~~\forall g\in G~$ and $~(h_{1}, h_{2}) \in H_{(2)}.$\\
$(3)~~$  For all $ h\in H~$  there exists an element $~\theta_{1}(g)~$  such that $~\theta_{1}(g)\odot h=h.$\\
$(4)~~g\odot \theta_{2}(h)=\theta_{2}(h),~$ for all  $~g\in G~$ and  $~h\in H.~$}
\end{defn}
\begin{Th} The semidirect product of two almost groupoids  is also an almost groupoid.
\end{Th}
{\it Proof.} Consider the almost groupoids  $ (G, \theta_{1}, \otimes_{1}, \iota_{1}, G_{0}) $ and  $ (H, \theta_{2}, \otimes_{2}, \iota_{2}, H_{0}).~$  

We will define an almost groupoid structure  on  the set $~ G \times H~$ over $~ G_{0} \times H_{0},~$  using an action of almost groupoid $~G~$ on the almost groupoid $~H.~$  For this purpose, let \\
$~\odot : G \times H \rightarrow H,~(g,h) \to g\odot h~$ be an almost groupoid action of $~G~$ on $~H.$\\
Let us now define the multiplication (denoted with "$\cdot$") on  $~ G\times H~$ as follows:\\
$(4.1)~~~(g_{1}, h_{1})\cdot (g_{2}, h_{2}):=(g_{1}\otimes_{1} g_{2}, h_{1} \otimes_{2}(g_{2}\odot h_{2})), ~\forall (g_{1}, g_{2})\in G_{(2)}~$ and $~(h_{1}, g_{2}\odot h_{2})\in H_{(2)}.$\\
According to the multiplication given in the relation $~(4.1)~,$ i the structure functions $~\theta: G\times H \rightarrow G_{0}\times H_{0}~$ and $~\iota: G\times H \rightarrow G\times H~$ are given by \\
$(4.2)~~~\theta (g,h):=(\theta_{1}(g), \theta_{2}(h))~~~\forall g\in  G~$ and $~~h\in H.$\\
$(4.3)~~~\iota (g,h):=(\iota_{1}(g), \iota_{1}(g)\odot (\iota_{2}(h)), \forall g\in G~$ and $~h\in H.$\\
Now, it is easy to verify that the conditions $ ({\bf AG1}) - {\bf AG3}) $ in Definition 4.1 hold.\\
First, we consider $~(g_{1}, h_{1}),  (g_{2}, h_{2}), (g_{3}, h_{3})\in G\times H.~$  We have:\\
$ (1)~~((g_{1}, h_{1})\cdot (g_{2}, h_{2}))\cdot (g_{3}, h_{3})= (g_{1}\otimes_{1} g_{2}, h_{1}\otimes_{2} (g_{1}\odot h_{2}))\cdot (g_{3}, h_{3})=$\\
$ (( g_{1}\otimes_{1} g_{2})\otimes_{1} g_{3}, h_{1}\otimes_{2} (g_{1}\odot h_{2}))\cdot (g_{3}, h_{3}) $=
$~ ( g_{1}\otimes_{1} g_{2}\otimes_{1} g_{3}, h_{1}\otimes_{2}(g_{1}\odot(h_{2}\otimes_{2}(g_{2}\odot h_{3})))~$  and\\
$ (g_{1}, h_{1})\cdot ((g_{2}, h_{2})\cdot (g_{3}, h_{3}))= (g_{1}, h_{1})\cdot (g_{2}\otimes_{1} g_{3}, h_{2}\otimes_{2} (g_{2}\odot h_{3}))=$\\
$  ( g_{1}\otimes_{1}( g_{2}\otimes_{1} g_{3}), h_{1}\otimes_{2} (g_{1}\odot (h_{2}\otimes_{2}(g_{2}\odot h_{3})))=~  ( g_{1}\otimes_{1} g_{2}\otimes_{1} g_{3},  h_{1}\otimes_{2} (g_{1}\odot (h_{2}\otimes_{2}(g_{2}\odot h_{3}))).$\\
Hence, $ ~((g_{1}, h_{1})\cdot (g_{2}, h_{2}))\cdot (g_{3}, h_{3})= (g_{1}, h_{1})\cdot ((g_{2}, h_{2})\cdot (g_{3}, h_{3})).$\\
$ (2)~~(g, h)\cdot \theta(g,h)=(g, h)\cdot (\theta_{1}(g), \theta_{2}(h))=(g\otimes_{1} \theta_{1}(g), h\otimes_{2} (g\odot \theta_{2}(h)))=$\\
$( g, h\otimes_{2} \theta_{2}(h))=(g,h)~$  and \\
$ \theta(g,h)\cdot (g, h)= (\theta_{1}(g), \theta_{2}(h))\cdot (g, h)=(\theta_{1}(g)\otimes_{1} g, \theta_{2}(h)\otimes_{2}( \theta_{1}(g)\odot h))=$\\
$( g, \theta_{2}(h)\otimes_{2} h)=(g,h).$\\
Hence, $ ~(g, h)\cdot \theta(g,h)= \theta(g,h)\cdot (g, h)=(g,h).$\\
$ (3)~~~(g, h)\cdot \iota(g,h)=(g, h)\cdot (\iota_{1}(g), (\iota_{1}(g)\odot \iota_{2}(h))=(g\otimes_{1} \iota_{1}(g), h\otimes_{2}(g\odot (\iota_{1}(g)\odot \iota_{2}(h)))=\\$
$(\theta_{1}(g), h\otimes_{2}(\theta_{1}(g)\odot \iota_{2}(h))= (\theta_{1}(g), h\otimes_{2}\iota_{2}(h))=(\theta_{1}(g), \theta_{2}(h))=\theta(g,h) ~$   and \\
$\iota(g,h)\cdot (g,h)=(\iota_{1}(g), \iota_{1}(g)\odot \iota_{2}(h))\cdot (g,h)=(\iota_{1}(g)\otimes_{1} g, (\iota_{1}(g)\odot \iota_{2}(h))\otimes_{2} (\iota_{1}(g)\odot h))=$\\
$(\theta_{1}(g), \iota_{1}(g)\odot (\iota_{2}(h)\otimes_{2} h))= (\theta_{1}(g), \iota_{1}(g)\odot \theta_{2}(h))=(\theta_{1}(g), \theta_{2}(h))=\theta(g,h).$\\
Hence    $~(g, h)\cdot \iota(g,h)=\iota(g,h)\cdot (g,h) = \theta(g,h).$\\ 
We conclude that $~ (G \times H, \theta, \cdot, \iota, G_{0}\times H_{0})~$  is an almost groupoid. \hfill$\Box$\\
The almost groupoid $~ G\times H~$  constructed in the proof of Theorem 4.1 is denoted by $~~G\times_{\propto} H~$ and is called the {\it semidirect product} of the almost groupoid $~G~$ with the almost groupoid $~H~$ via the almost groupoid action 
$~\odot: G\times H \rightarrow H.$\\[-0.7cm]

\vspace*{0.2cm}

Author's adress\\[0.2cm]
\hspace*{0.7cm}West University of Timi\c soara. Seminarul de Geometrie \c si Topologie\\
\hspace*{0.7cm} Department of Mathematics\\
\hspace*{0.7cm} Bd. V. P{\^a}rvan,no.4, 300223, Timi\c soara, Romania\\
\hspace*{0.7cm}E-mail: gheorghe.ivan@e-uvt.ro\\

\end{document}